\documentclass[12pt]{article}%
\usepackage{amsmath}
\usepackage{graphicx}
\usepackage{amsfonts}
\usepackage{amssymb}
\usepackage{latexsym}
\usepackage{color}
\usepackage[colorlinks,pdfpagemode=FullScreen,pdftex,debug]{hyperref}
\usepackage[T1]{fontenc}
\usepackage[sc]{mathpazo}
\usepackage{microtype}%
\newtheorem{theorem}{Theorem}

\newtheorem{conjecture}[theorem]{Conjecture}

\newtheorem{problem}[theorem]{Problem}
\newtheorem{proposition}[theorem]{Proposition}

\newenvironment{proof}[1][Proof]{\noindent{\textbf {#1}  }}  {\hfill$\Box$\bigskip}
\makeatletter
\definecolor {ttlcol}{rgb}{.96,0,.12}
\definecolor {ttlcol}{rgb}{.85,0,0}
\definecolor {pagecol}{rgb}{.98,.98,1}
\definecolor {thmcol}{rgb}{.6,0,0}
\definecolor {deficol}{rgb}{0,.5,.25}
\definecolor {quecol}{rgb}{0,.25,.5}
\definecolor {bgnd}{rgb}{1,1,.95}
\begin{document}

\title{\textbf{The energy of }$C_{4}$\textbf{-free graphs of bounded degree}}
\author{V. Nikiforov\\{\small Department of Mathematical Sciences, University of Memphis, }\\{\small Memphis TN 38152, USA, e-mail: vnikifrv@memphis.edu}}
\maketitle

\begin{abstract}
Answering some questions of Gutman, we show that, except for four specific
trees, every connected graph $G$ of order $n,$ with no cycle of order $4$ and
with maximum degree at most $3,$ satisfies
\[
\left\vert \mu_{1}\right\vert +\cdots+\left\vert \mu_{n}\right\vert \geq n,
\]
where $\mu_{1},\ldots,\mu_{n}$ are the eigenvalues of $G$.

We give some general results and state two conjectures.\bigskip



\end{abstract}

Our notation follows \cite{Bol98} and \cite{CDS80}; in particular, we write
$V\left(  G\right)  $ for the vertex set of a graph $G$ and $\left\vert
G\right\vert $ for $\left\vert V\left(  G\right)  \right\vert .$ Also,
$e\left(  G\right)  $ stands for the number of edges of $G,$ and
$\Delta\left(  G\right)  $ for its maximum degree.

Let $G$ be a graph on $n$ vertices and $\mu_{1}\geq\cdots\geq\mu_{n}$ be the
eigenvalues of its adjacency matrix. The value $\mathcal{E}\left(  G\right)
=\left\vert \mu_{1}\right\vert +\cdots+\left\vert \mu_{n}\right\vert ,$ called
the \emph{energy }of $G,$ has been studied intensively - see \cite{Gut05} for
a survey.

Motivated by questions in theoretical chemistry, Gutman \cite{Gut07} initiated
the study of connected graphs satisfying\emph{ }$\mathcal{E}\left(  G\right)
\geq\left\vert G\right\vert ,$ raising some problems, whose simplest versions
read as:

\begin{problem}
Characterize all trees $T$ with $\Delta\left(  T\right)  \leq3$ satisfying
$\mathcal{E}\left(  T\right)  <\left\vert T\right\vert .$
\end{problem}

\begin{problem}
Characterize all connected graphs $G$ with $\Delta\left(  G\right)  \leq3$
satisfying $\mathcal{E}\left(  G\right)  <\left\vert G\right\vert .$
\end{problem}

Here we show that if $G$ is a connected $C_{4}$-free graph such that
$\Delta\left(  G\right)  \leq3$ and $\mathcal{E}\left(  G\right)  <\left\vert
G\right\vert ,$ then $G$ is one of four exceptional trees. This completely
solves the first problem and partially the second one.

Let $d\geq3$ and $\alpha\left(  d\right)  $ be the largest root of the
equation%
\[
4x^{3}-\left(  2d+1\right)  x+d=0\text{.}%
\]

\begin{theorem}
\label{th1}Let $G$ be a $C_{4}$-free graph with no isolated vertices. If
$e\left(  G\right)  \geq\alpha\left(  d\right)  \left\vert G\right\vert $ and
$\Delta\left(  G\right)  \leq d,$ then $\mathcal{E}\left(  G\right)
>\left\vert G\right\vert .$
\end{theorem}

We note first that Theorem \ref{th1} implies Theorem 1 of \cite{Gut07}, but
the check of this implication is somewhat involved.

To prove Theorem \ref{th1} we need three propositions. The first one is known
and its proof is omitted.

\begin{proposition}
\label{pro1}Let $G$ be a graph of order $n,$ $C$ be the number of its
$4$-cycles, and $\mu_{1},\ldots,\mu_{n}$ be its eigenvalues. Then%
\begin{align*}
\mu_{1}^{2}+\cdots+\mu_{n}^{2}  &  =2e\left(  G\right)  ,\\
\mu_{1}^{4}+\cdots+\mu_{n}^{4}  &  =2%
%TCIMACRO{\tsum \limits_{u\in V\left(  G\right)  }}%
%BeginExpansion
{\textstyle\sum\limits_{u\in V\left(  G\right)  }}
%EndExpansion
d_{u}^{2}-2e\left(  G\right)  +8C.
\end{align*}

\end{proposition}

Next we give a simple bound on the sum of squares of degrees in graphs.

\begin{proposition}
\label{pro2}Let $G$ be a graph with $n$ vertices and $m$ edges, with no
isolated vertices, and let $d_{1},\ldots,d_{n}$ be its degrees. If
$\Delta\left(  G\right)  \leq d,$ then%
\[
d_{1}^{2}+\cdots+d_{n}^{2}\leq\left(  2d+2\right)  m-dn.
\]

\end{proposition}

\begin{proof}
Summing the inequality $\left(  d_{i}-1\right)  \left(  d_{i}-d\right)  \leq0$
for $i=1,\ldots,n,$ we find that
\[
d_{1}^{2}+\cdots+d_{n}^{2}-d_{1}-\cdots-d_{n}-d\left(  d_{1}+\cdots
+d_{n}\right)  +dn\leq0,
\]
and the desired inequality follows in view of $d_{1}+\cdots+d_{n}=2m$.
\end{proof}

The following proposition gives more explicit relations between $d$ and
$\alpha\left(  d\right)  $.

\begin{proposition}
If $d=3,$ then $\alpha\left(  d\right)  =1.$ If $d\geq4,$ then
\begin{equation}
\sqrt{\left(  2d+1\right)  /4}-1/3<\alpha\left(  d\right)  <\sqrt{\left(
2d+1\right)  /4}. \label{abnds}%
\end{equation}

\end{proposition}

\begin{proof}
If $d=3,$ we have%
\[
4x^{3}-7x+3=4x\left(  x-1\right)  \left(  x+1\right)  -3\left(  x-1\right)
=\left(  x-1\right)  \left(  2x-1\right)  \left(  2x+3\right)
\]
and the first assertion follows.

If $x\geq\sqrt{\left(  2d+1\right)  /4},$ we have
\[
4x^{3}-\left(  2d+1\right)  x+d\geq\left(  \left(  2d+1\right)  -\left(
2d+1\right)  \right)  x+d>0,
\]
so the upper bound in (\ref{abnds}) follows. On the other hand,%
\[
4\left(  \sqrt{\frac{2d+1}{4}}-\frac{1}{3}\right)  ^{3}-\left(  2d+1\right)
\left(  \sqrt{\frac{2d+1}{4}}-\frac{1}{3}\right)  +d<-\frac{19d+11}{3}%
+\frac{1}{3}\sqrt{8d+4}<0,
\]
implying the lower bound in (\ref{abnds}) and completing the proof.
\end{proof}

\bigskip

\begin{proof}
[\textbf{Proof of Theorem \ref{th1}}]As noted by Rada and Tineo \cite{RaTi04}%
,
\begin{equation}
\left(  \left\vert \mu_{1}\right\vert +\cdots+\left\vert \mu_{n}\right\vert
\right)  ^{2/3}\left(  \mu_{1}^{4}+\cdots+\mu_{n}^{4}\right)  ^{1/3}\geq
\mu_{1}^{2}+\cdots+\mu_{n}^{2}. \label{RTin}%
\end{equation}
We first show that, in our case, inequality (\ref{RTin}) is strict. Indeed, it
is a particular case of H\"{o}lder's inequality; hence, if equality holds, the
vectors $\left(  \left\vert \mu_{1}\right\vert ,\ldots,\left\vert \mu
_{n}\right\vert \right)  $ and $\left(  \mu_{1}^{4},\ldots,\mu_{n}^{4}\right)
$ are linearly dependent. That is to say, all nonzero eigenvalues of $G$ have
the same absolute value, and so, $G$ is a union of complete bipartite graphs.
Since $G$ is $C_{4}$-free, $G$ must be a forest, contradicting the premise
$e\left(  G\right)  \geq\alpha\left(  d\right)  n\geq n,$ and proving the claim.

From (\ref{RTin}), by Propositions \ref{pro1} and \ref{pro2}, we obtain
\[
\mathcal{E}^{2}\left(  G\right)  >\frac{\mu_{1}^{2}+\cdots+\mu_{n}^{2}}%
{\mu_{1}^{4}+\cdots+\mu_{n}^{4}}=\frac{8m^{3}}{2%
%TCIMACRO{\tsum \limits_{u\in V\left(  G\right)  }}%
%BeginExpansion
{\textstyle\sum\limits_{u\in V\left(  G\right)  }}
%EndExpansion
d^{2}\left(  u\right)  -2m}\geq\frac{4m^{3}}{\left(  2d+1\right)  m-dn}%
=n^{2}\frac{4\left(  m/n\right)  ^{3}}{\left(  2d+1\right)  \left(
m/n\right)  -d}.
\]
Using calculus, we find that the expression
\[
\frac{4x^{3}}{\left(  2d+1\right)  x-d}%
\]
increases with $x.$ Hence, the premise $m/n\geq\alpha\left(  d\right)  $
implies that%
\[
\mathcal{E}^{2}\left(  G\right)  >n^{2}\frac{4\alpha^{3}}{\left(  2d+1\right)
\alpha-d}=n^{2},
\]
completing the proof.
\end{proof}

\bigskip

As a corollary we obtain the following

\begin{theorem}
\label{th2}Let $G$ be a graph of order $n$ with at least $n$ edges and with no
isolated vertices. If $G$ is $C_{4}$-free and $\Delta\left(  G\right)  \leq3,$
then $\mathcal{E}\left(  G\right)  >n.$
\end{theorem}

It is reasonable to conjecture that $\mathcal{E}\left(  G\right)
\geq\left\vert G\right\vert $ holds for all connected $C_{4}$-free graphs with
$\Delta\left(  G\right)  \leq3.$ In view of Theorem \ref{th1}, this assertion
can fail only if $G$ is a tree, and indeed, as pointed by Gutman, there are
four trees for which the assertion fails.

\begin{fact}
[\textbf{Gutman \cite{Gut07}}]\label{f1}The trees $K_{1},$ $K_{1,2},$
$K_{1,3}$, and the balanced binary tree on $7$ vertices are the only trees of
order up to $22$ with $\Delta\left(  T\right)  \leq3$ and $\mathcal{E}\left(
T\right)  <\left\vert T\right\vert $.$\hfill\square$
\end{fact}

However, it turns out that these trees are the only exceptions, and the
following theorem holds.

\begin{theorem}
\label{th3}Let $T$ be a tree different from the four trees listed in Fact
\ref{f1}$.$ If $\Delta\left(  T\right)  \leq3,$ then $\mathcal{E}\left(
T\right)  \geq\left\vert T\right\vert $.
\end{theorem}

\begin{proof}
In view of Fact \ref{f1}, we shall assume that $n\geq23$. Let $\mu_{1}%
,\ldots,\mu_{n}$ be the eigenvalues of $T.$ Set $\mu=\mu_{1}$ and note that,
since $T$ is bipartite, we also have $\left\vert \mu_{n}\right\vert =\mu.$
H\"{o}lder's inequality implies that
\[
\left(  \left\vert \mu_{2}\right\vert +\cdots+\left\vert \mu_{n-1}\right\vert
\right)  ^{2/3}\left(  \mu_{2}^{4}+\cdots+\mu_{n-1}^{4}\right)  ^{1/3}\geq
\mu_{2}^{2}+\cdots+\mu_{n-1}^{2}.
\]
Hence Propositions \ref{pro1} and \ref{pro2} give
\begin{equation}
\left(  \mathcal{E}\left(  G\right)  -2\mu\right)  ^{2}\geq\frac{\mu_{2}%
^{2}+\cdots+\mu_{n-1}^{2}}{\mu_{2}^{4}+\cdots+\mu_{n-1}^{4}}=\frac{\left(
2m-2\mu^{2}\right)  ^{3}}{2%
%TCIMACRO{\tsum \limits_{u\in V\left(  G\right)  }}%
%BeginExpansion
{\textstyle\sum\limits_{u\in V\left(  G\right)  }}
%EndExpansion
d^{2}\left(  u\right)  -2m-2\mu^{4}}\geq\frac{4\left(  n-1-\mu^{2}\right)
^{3}}{4n-7-\mu^{4}}. \label{eq1}%
\end{equation}
First we show that if $\mu\geq\sqrt{7}$, then
\begin{equation}
\frac{4\left(  n-\left(  \mu^{2}+1\right)  \right)  ^{3}}{4n-7-\mu^{4}%
}>\left(  n-2\mu\right)  ^{2}. \label{eq2}%
\end{equation}
We shall make use of the fact $\mu<\Delta\left(  T\right)  \leq3.$ After some
algebra we see that (\ref{eq2}) is equivalent to%
\[
\left(  \mu^{4}-12\mu^{2}+16\mu-5\right)  n^{2}+4\left(  -\mu^{5}+3\mu
^{4}+2\mu^{2}-7\mu+3\right)  n+4\left(  -3\mu^{4}+4\mu^{2}-1\right)  >0.
\]
We have, in view of $\sqrt{7}\leq\mu\leq3,$
\begin{align*}
&  \left(  \mu^{4}-12\mu^{2}+16\mu-5\right)  n^{2}+4\left(  -\mu^{5}+3\mu
^{4}+2\mu^{2}-7\mu+3\right)  n+4\left(  -3\mu^{4}+4\mu^{2}-1\right) \\
&  >\left(  7\left(  7-12\right)  +16\sqrt{7}-5\right)  n^{2}+4\left(
2\mu^{2}-7\mu+3\right)  n+4\left(  -243+27\right) \\
&  >2n^{2}+4\left(  14-7\sqrt{7}+3\right)  n-4\cdot216\geq2n^{2}%
-8n-4\cdot216\\
&  \geq46\cdot19-4\cdot216=10>0.
\end{align*}
Combining (\ref{eq1}) and (\ref{eq2}), we complete the proof if $\mu\geq
\sqrt{7}.$

Assume now that $\mu<\sqrt{7}.$ Hofmeister's inequality implies that
\[%
%TCIMACRO{\tsum \limits_{u\in V\left(  G\right)  }}%
%BeginExpansion
{\textstyle\sum\limits_{u\in V\left(  G\right)  }}
%EndExpansion
d^{2}\left(  u\right)  \leq n\mu^{2}<7n,
\]
and, as in the proof of Theorem \ref{th1}, we obtain%
\[
\mathcal{E}^{2}\left(  G\right)  \geq\frac{8\left(  n-1\right)  ^{3}}{2%
%TCIMACRO{\tsum \limits_{u\in V\left(  G\right)  }}%
%BeginExpansion
{\textstyle\sum\limits_{u\in V\left(  G\right)  }}
%EndExpansion
d^{2}\left(  u\right)  -2\left(  n-1\right)  }\geq\frac{8\left(  n-1\right)
^{3}}{14n-2\left(  n-1\right)  }=\frac{8\left(  n-1\right)  ^{3}}{12n+2}%
>n^{2},
\]
completing the proof.
\end{proof}

\subsubsection*{Concluding remarks}

Some of the above results can be strengthened. Here we formulate two plausible conjectures.

Let the tree $\mathcal{B}_{n}$ be constructed by taking three disjoint copies
of the balanced binary tree of order $2^{n+1}-1$ and joining an additional
vertex to their roots.

\begin{conjecture}
The limit
\[
c=\lim_{n\rightarrow\infty}\frac{\mathcal{E}\left(  \mathcal{B}_{n}\right)
}{3\cdot2^{n+1}-2}%
\]
exists and $c>1.$
\end{conjecture}

\begin{conjecture}
Let $\varepsilon>0.$ If $T$ is a sufficiently large tree with $\Delta\left(
T\right)  \leq3$, then $\mathcal{E}\left(  T\right)  \geq\left(
c-\varepsilon\right)  \left\vert T\right\vert .$
\end{conjecture}

The empirical data given in \cite{Gut07} seem to corroborate these
conjectures, but apparently new techniques are necessary to prove or disprove
them.

\end{document}